\documentclass[11pt]{article}
\usepackage{latexsym}
\usepackage{amsfonts}

\title{An historical overview of the influence of technology on mathematical competitions \\[.4in]}

\author{B\'{e}la Bajnok \\[.1in] Department of Mathematics \\ Gettysburg College \\
Gettysburg, PA 17325-1486 USA \\E-mail:  bbajnok@gettysburg.edu }

\date{April 19, 2010}

\begin{document}

\maketitle

\thispagestyle{empty}

\begin{abstract}

We provide an historical overview of how advances in technology influenced high school and university mathematical competitions in the United States and at the International Mathematical Olympiad.  While students are not allowed the usage of technological aids during mathematical competitions, the developments in technology (especially graphing technology) throughout the past century and the increasing employment of such aids in the classroom have affected both the nature of the proposed problems and their expected solutions.  We examine several interesting examples from competitions going back several decades.

\end{abstract}

\section{Introduction}

Our goal in this paper is to investigate how advances in technology influenced mathematical problem posing and problem solving.  We will concentrate on high school and university mathematical competitions in the United States and at the International Mathematical Olympiad.  The 10 problems we chose demonstrate an increasing level of technological sophistication when it comes to both problem posing and problem solving.   

\section{About mathematical competitions available in the United States}

Let us first briefly review mathematical contests in the US.

The American Mathematics Competitions (AMC) has a long and distinguished history (although not nearly as long and as distinguished as similar competitions in other parts of the world, especially in Eastern and Central Europe).  The first American mathematical contest took place in 1950, sponsored by the New York Metropolitan Section of the Mathematics Association of America (MAA).  

Today, the AMC organizes three competitions: the American Mathematics Contest (AMC 8, 10, 12, for grades 8, 10, and 12, respectively), the American Invitational Mathematics Examination (AIME), and the United States of America Mathematical Olympiad (USAMO).  Each year over 400,000 students in over 5,000 schools participate in the AMC Contests. Of these, about 10,000 students qualify to participate in the AIME. From this group approximately 500 students are invited to take the prestigious USAMO. 

The current three-level AMC system replaced a previous two-tier system (AJHSME - American Junior High School Mathematics Examination and AHSME - American High School Mathematics Examination) in 2000.  Here we focus on AMC 12 (and its predecessor, the AHSME), which is currently a 25-question, 75-minute multiple choice examination in secondary school mathematics containing problems which can be understood and solved with pre-calculus concepts.  There is no penalty for guessing. 

The American Invitational Mathematics Examination (AIME), is a 15 -question, 3-hour examination in which each answer is an integer number from 0 to 999.  The questions on the AIME are much more difficult than those at the AMC, and students are very unlikely to obtain the correct answer by guessing.

The United States of America Mathematical Olympiad (USAMO), is a 6-question, 2-day, 9-hour essay/proof examination. The USAMO problems are meant to be very challenging but, like other AMC competition problems, all USAMO problems can be solved with pre-calculus methods.

The AMC culminates with the International Mathematical Olympiad (IMO) which is a 10-14 day trip and contest for the top 6 high school students at the Mathematical Olympiad Summer Program (MOSP), who comprise the United States IMO team and represent the USA at the IMO.  The IMO began in 1959 in Romania for countries of the Warsaw Pact; today there are nearly 100 countries represented.  The USA has participated since 1974.

Finally, the mathematical competition for college students: the William Lowell Putnam Mathematical Competition or Putnam Exam.  This is an annual mathematics competition for undergraduate college students of the United States and Canada, awarding scholarships and cash prizes.  The Putnam Exam is a 6-hour, 12-question competition, which was first offered in 1938.

\section{The beginning: the challenges of graphing in the 1950's}

Before the rapid spread of graphing and other technological aids in the classroom, the American curriculum considered even trivial graphing as a challenge.  It is therefore not surprising that very few competition problems involved graphs or graphical approaches.  While published solutions occasionally refer to graphs, they almost never show figures.  Let us look at a few examples from this age.

Problem 44 on the 1950 AHSME contest asks whether the graph of $y=\log x$ ``cuts'' the $x$-axis, ``cuts'' the $y$-axis, ``cuts'' all lines parallel to the $x$-axis, etc.  We may be surprised at both the level of the question (Problem 44 was considered one of the most challenging questions at a time when the contest contained 50 questions) and the fact that the published solution (\cite{Sal:1961a}) does not refer to graphs.

Four years later, Problem 17 on the 1954 AHSME exam asks whether the function $f(x)=2x^3-7$ ``(A) goes up to the right and down to the left'' or ``(B) down to the right and up to the left,'' etc.  The official solution (\cite{Sal:1961a}) is the following:  ``There are a number of ways determining that (A) is the correct answer.  One (obvious) way is to plot a sufficient number of points.  Note: Graphs of cubics, generally, are beyond the indicated scope.'' Of course, no graph is given in the solution.

Problem 42 of the same exam is especially noteworthy in our discussion from the problem posing point of view.  It asks students to consider the graphs of (1) $y=x^2-0.5x+2$ and (2) $y=x^2+0.5x+2$ in the same coordinate system.  The problem states that ``these parabolas have exactly the same shape'' (\cite{Sal:1961a}), which may surprise students who correctly recall that all parabolas are similar anyway.  The task is then to decide if ``(A) the graphs coincide; (B) the graph of (1) is lower than the graph of (2); (C) the graph of (1) is to the left of the graph of (2);  (C) the graph of (1) is to the right of the graph of (2); or (E) the graph of (1) is higher than the graph of (2).''  While we clearly see that the each parabola is a horizontal shift from the other, we find this question very poorly phrased.  (Fortunately, there were no questions to decide whether $y=x-1$ is below or to the right of $y=x+1$ or if the sine curve is to the left or to the right of the cosine curve.)

\section{Increased expectations in the 1960's}

The level of difficulty of high school competitions was raised greatly during the 1960's.  Official solutions, however, still lacked graphing approaches.  Let us look at some examples from this period.

Problem 33 of the 1966 AHSME competition asked students to decide how many distinct solutions the equation
$$\frac{x-a}{b}+\frac{x-b}{a}=\frac{b}{x-a}+\frac{a}{x-b}$$ had for $x$.  The official solution (\cite{SalEar:1973a}) suggests an algebraic solution involving case separation.  We can simply realize, however, that the left-hand side of the equation is a linear function of $x$, while the right-hand side is a rational function with two vertical and one horizontal asymptotes ($x=a$, $x=b$, and $y=0$).  Further analyzing these asymptotes, we can quickly conclude that the equation has three distinct solutions.  (As this is a question on a multiple choice exam where the answer does not depend on parameters, we can freely choose our favorite $a$ and $b$ to investigate the graphs.)

The second problem on the 1968 IMO asked contestants to find all natural numbers $x$ for which the product of the digits of $x$ equals $x^2-10x -22$.  The original solution (\cite{Gre:1978a}) notes that the product of the digits of $x$ is at most $x$, then, after some algebraic manipulations and case separations, concludes that the only solution is $x=12$.  We can, however, note that if $x^2-10x-22$ is the product of the digits of $x$ then it has to be nonnegative, and analyze the behavior of the corresponding parabola to find that it is between 0 and $x$ only when $x=12$.  We can check that $x=12$ indeed works. 

\section{Technological advances of the 1970's and 1980's}

Graphing technology and graphical approaches became integral parts of the American high school curriculum during the late 1970's and early 1980's.  Problem posers of competitions started to reinforce graphing, and official solutions included graphical approaches.

Let us mention here two examples only.  Problem 13 on the 1979 AHSME competition asks for the solution set of the inequality $$y-x < \sqrt{x^2},$$ and Problem 18 on the 1981 AHSME contest asks for the number of real solutions to the equation $$\frac{x}{100}=\sin x.$$
The solution to both problems (\cite{Art:1985a}) are graphical in nature and are accompanied by figures.

\section{Graphing technology matures during the 1990's} 

During the 1990's, graphs and graphical techniques are becoming to be recognized as frequently more efficient than routine algorithms even when routine algorithms are available.

Consider, for example, Problem 26 on the 1993 AHSME contest.  The question asks for the maximum value of $$f(x)=\sqrt{8x-x^2}-\sqrt{14x-x^2-48}.$$
The two published solutions in \cite{AHS} both note that $x$ has to lie in the interval $[6.8]$.  While a second, purely algebraic, solution is also given where the technique of multiplication by the conjugate is used, the first solution notes that the graphs corresponding to the two terms in $f(x)$ are semi-circles and thereby finds that the difference is maximal when $x=6$.  Another potential solution (not given) would simply analyze the graphs of the two parabolas $y=8x-x^2$ and $y=14x-x^2-48$.

For another example, let us mention Problem A2 on the 1994 Putnam Exam: ``Let $A$ be the area of the region in the first quadrant bounded by the line $y=x/2$, the $x$-axis, and the ellipse $$\frac{x^2}{9}+y^2=1.$$  Find the positive number $m$ such that $A$ is equal to the area of the region in the first quadrant bounded by the line $y=mx$, the $y$-axis, and the [same] ellipse.''  The only reasonable solution to this problem (\cite{Klo:1995a}) is to apply a linear transformation which takes the ellipse to the unit circle (and the lines $y=x/2$ and $y=mx$ to the lines $y=3x/2$ and $y=3mx$, respectively).  The value of $m$ can now be easily found using symmetry of the circle about the line $y=x$. While we do not have data on how many contestants attempted to solve the problems this way rather than trying to use integral calculus, we expect that many of them did, particularly students who are used to viewing the circle on a calculator or computer screen as an ellipse due to having scales of different lengths on the two axes.

\section{The ubiquity of technological tools  today} 

These days, it is rare to meet students at the high school or college level who are not familiar with technological tools in the classroom.  In fact, in this author's experience, their reliance on graphical and numerical sources almost overwhelm them; they often reach for these tools before even fully understanding the questions.  Nevertheless, this has had a very positive affect on the expectations problem posers have for contestants as far as graphical experiences.  We demonstrate this with an example from this year.  

Problem 21 of 2010 AMC 12 is as follows. ``The graph of $y=x^6-10x^5+29x^4-4x^3+ax^2$ lies above the line $y=bx+c$ except at three values of $x$, where the graph and the line intersect. What is the largest of these values?''

The solution requires students to understand that the function $$f(x)=x^6-10x^5+29x^4-4x^3+ax^2-bx-c$$ lies entirely on or above the $x$-axis and has three double roots.  Therefore, $f(x)$ is the square of the cubic polynomial $$g(x)=(x-p)(x-q)(x-r).$$  A bit of algebra then easily yields that $\{p,q,r\}=\{-1,2,4\}$ and thus the correct answer is $x=4$.

Here both the question and the solution is at a level that is substantially beyond what one would have expected just a few years ago.

\section{Conclusions}

As our ``top ten'' list of exemplary problems demonstrates, technological advances have had a great influence on both mathematical problem posing and problem solving in American education.  Most of the graphing approaches we found in our research are not dependent on specific software tools or calculator features.  Instead, students whose mathematical experiences included technological tools found a way to use these skills and knowledge in creative work, whether technological assistance is available for them at the time or not.  This demonstrates that with modern tools students at any level can finally reach a better understanding of the subject.

In closing, we pose a few questions.  Now that several generations of teachers and students have been trained to use technological tools, how will the future of mathematical competitions be reshaped?  Will certain areas of mathematics receive increased or decreased representation at contexts?  To what extent will the awareness of graphical approaches help explore different solutions to problems?  Is there such a thing as relying on technology ``too much''?  Will the nature of mathematical competitions change due to general availability of technological aids?  We can hardly wait to find out the answers.

{\bf Acknowledgments.}  A previous version of this paper, by the current author and \'Agnes Tuska (\cite{BajTus:1996a}), appeared in {\em Proceedings of the Eighth Annual International Conference on Technology in Collegiate Mathematics}.  The author expresses his gratitude to \'Agi for the many interesting discussions on this subject.

\end{document}